\documentclass[a4paper, 12pt]{article}

\usepackage{geometry}                
\usepackage{calligra}
\usepackage[T1]{fontenc}
\usepackage[english]{babel}
\usepackage{amsmath}
\usepackage{amsthm}
\usepackage{amssymb}
\usepackage{enumerate}
\usepackage{xspace}
\usepackage{euscript}
\usepackage{graphicx}
\usepackage{amscd}
\usepackage{epsfig}
\usepackage{tabularx}
\usepackage{pgfplots}
\usepackage{tikz}
\usetikzlibrary{datavisualization}
\usetikzlibrary{datavisualization.formats.functions}

\newtheorem{theorem}{Theorem}
\newtheorem{remark}{Remark}
\newtheorem{definition}{Definition}
\newtheorem{proposition}{Proposition}
\newcommand{\supp}{{\rm Supp} \,}
\newcommand{\N}{{\mathbb N}}

\newcommand{\R}{{\mathbb R}}

\title{On backward uniqueness for parabolic equations when Osgood continuity of the coefficients fails at one point}
\author{Daniele Del Santo \and Martino Prizzi}
\date{\null}                                           

\begin{document}
\maketitle

\begin{abstract}
We prove the uniqueness for backward parabolic equations whose coefficients are Osgood continuous in time for $t>0$ but not at  $t=0$.
\end{abstract}




\section{Introduction}

Let us consider the following backward parabolic operator
$$
L=\partial_t+ \sum_{i, j = 1}^n \partial_{x_j} (a_{jk} (t,x)
 \partial_{x_k} )+\sum_{ j = 1}^n b_j(t,x)\partial_{x_j}+c(t,x).
$$
We assume that all coefficients are  defined in $[0,T]
\times {\mathbb R}^n_x$,
measurable  and bounded; $( a_{jk}(t,x))_{jk}$ is a real symmetric matrix for
all $(t,x)\in [0,T]\times
{\mathbb R}^n_x$ and there exists $\lambda_0\in (0,1]$ such that
$$
\sum_{j, k = 1}^n  a_{jk} (t,x)
\xi_j \xi_k\geq \lambda_0|\xi|^2\label{k2}
$$ for all $(t,x)\in [0,T]\times  {\mathbb R}^n_x$ and $\xi\in {\mathbb R}^n_\xi$.

Given a functional space ${\cal H}$ we say that the operator $L$ has the ${\cal
H}$--uniqueness property if, whenever $u\in
{\cal H}$, $Lu=0$ in $[0,T]\times  {\mathbb R}^n_x$ and
$u(0,x)=0$ in ${\mathbb R}^n_x$, then $u=0$ in $[0,T]\times  {\mathbb R}^n_x$.
Our choice for ${\mathcal H}$ is the space of functions 
$$
 {\mathcal H}= H^1((0,T), L^2({\mathbb R}^n_x))\cap L^2( (0,T), H^2({\mathbb R}^n_x)).
 $$
This choice is natural, since 
it follows from elliptic regularity results that the domain of the operator $-\sum_{j,k = 1}^n \partial_{x_j} (a_{jk} (t,x)
\partial_{x_k} )$ in $L^2(\mathbb R^n)$ is $H^2(\mathbb R^n)$ for all $t\in[0,T]$.

In our previous papers \cite{DSP1,DSP2} we investigated the problem of finding the minimal
regularity assumptions on the coefficients
$a_{jk}$ ensuring the ${\cal H}$--uniqueness property to $L$. Namely, we proved
 the $ {\cal
H}$--uniqueness property for
the operator $L$
when the coefficients $a_{jk}$ are Lipschitz continuous in $x$
and the regularity in $t$ is given in
terms of a modulus of
continuity $\mu$, i. e.  
$$
\sup_{s_1,\,s_2\in[0,T],\atop x\in\R^n}\frac{|a_{j,k}(s_1,x)-a_{j,k}(s_2,x)|}{\mu(|s_1-s_2|)}\leq C,
$$
where $\mu$ satisfies the so called Osgood condition
$$
\int_0^1{1\over \mu(s)}\, ds=+\infty.
$$

Suitable counterexamples show that Osgood condition is sharp for backward uniqueness in parabolic equations: given any non-Osgood modulus of continuity $\mu$, it is possible to construct a backward parabolic equation, whose coefficients are $C^\infty$ in $x$ and $\mu$-continuous in $t$, for which the ${\cal H}$--uniqueness property does not hold. In the mentioned counterexamples the coefficients are in fact $C^\infty$ in $t$ for $t\not=0$, and 
Osgood continuity fails only at $t=0$.

\bigskip
In this paper we show that if the loss of Osgood continuity is properly controlled as $t\to 0$, then we can recover the ${\cal H}$--uniqueness property for $L$. Our hypothesis reads as follows: given a modulus of continuity $\mu$ satisfying the Osgood condition,
we assume that the coefficients $a_{jk}$ are H\"older continuous with respect to $t$ on $[0,T]$, and for all $t\in]0,T]$
\begin{equation}
\label{control}\sup_{s_1,\,s_2\in[t,T],\atop x\in\R^n}\frac{|a_{j,k}(s_1,x)-a_{j,k}(s_2,x)|}{\mu(|s_1-s_2|)}\leq C t^{-\beta},
\end{equation} 
where $0<\beta<1$. The cofficients $a_{jk}$ are assumed to be globally Lipschitz continuous in $x$. Under such hypothesis we prove that the ${\cal H}$--uniqueness property holds for $L$. As in our previous papers \cite{DSP1,DSP2}, the uniqueness result is consequence of a Carleman estimate with a weight function shaped on the modulus of continuity $\mu$.  The weight function is obtained as solution of a specific second order ordinary differential equation. In the previous results cited above, the corresponding o.d.e. is autonomous. Here, on the contrary, the time dependent control (\ref{control}) yields to a non-autonomous o.d.e.. Also, the "Osgood singularity" of $a_{jk}$ at $t=0$ introduces a number of new technical difficulties
 which are not present in the fully Osgood-regular situation considered before.

The result is sharp in the following sense: we exhibit a counterexample in which the coefficients $a_{jk}$ are H\"older continuous with respect to $t$ on $[0,T]$, for all $t\in]0,T]$ and for all $\epsilon>0$
\begin{equation}
\label{control2}\sup_{s_1,\,s_2\in[t,T],\atop x\in\R^n}\frac{|a_{j,k}(s_1,x)-a_{j,k}(s_2,x)|}{|s_1-s_2|}\leq C t^{-(1+\epsilon)},
\end{equation} 
and the operator $L$ does not have the ${\cal H}$--uniqueness property. The 
borderline case $\epsilon=0$ in (\ref{control2}) is considered in paper \cite{DJ}. In such a situation only a very particular uniqueness result holds and the problem remains essentially open. 



\section{Main result}

We start with the definition of modulus of continuity.

\begin{definition}
A function $\mu: [0,\, 1]\to [0,\,1]$ is a {\rm modulus of continuity} if it is continuous, concave, strictly increasing and $\mu(0)=0$, $\mu(1)=1$. 
\end{definition}

\begin{remark}
Let $\mu$ be a modulus of continuity. Then
\begin{itemize}
\item for all $s\in [0,\, 1]$, $\mu(s)\geq s$;
\item on $(0,\, 1]$, the function $s\mapsto \frac{\mu(s)}{s}$ is decreasing;
\item the limit $\lim_{s\to 0^+}  \frac{\mu(s)}{s}$ exists;
\item on $[1,\, +\infty[$, the function $\sigma\mapsto \sigma \mu (\frac{1}{\sigma})$ is increasing;
\item on $[1,\, +\infty[$, the function $\sigma\mapsto \frac{1}{\sigma^2 \mu(\frac{1}{\sigma})}$ is decreasing.
\end{itemize}
\end{remark}

\begin{definition}
Let $\mu$ be a modulus of continuity and let $\varphi:I \to B$, where $I$ is an interval in $\R$ and $B$ is a Banach space.
$\varphi$ is a {\rm function in $C^\mu(I,B)$} if $\varphi\in L^\infty(I,B)$ and
$$
\|\varphi\|_{C^\mu(I,B)}= \|\varphi\|_{L^\infty(I,B)}+ \sup_{t,s\in I\atop 0<|t-s|<1} \frac{\| \varphi(t)-\varphi(s)\|_B}{\mu(|t-s|)}<+\infty.
$$
\end{definition}
\begin{remark}
Let $\alpha\in(0,1)$ and $\mu(s)=s^\alpha $. Then $C^\mu(I,B)=C^{0,\alpha}(I,B)$, the space of H\"older-continuous functions. 
Let $\mu(s)=s $. Then  $C^\mu(I,B)=Lip (I,B)$, the space of bounded Lipschitz-continuous functions.
\end{remark}

We introduce the notion of Osgood modulus of continuity.

\begin{definition}
Let $\mu$ be a modulus of continuity. $\mu$ {\rm satisfies the Osgood condition} if 
\begin{equation}\label{osgood}
\int_0^1 \frac{1}{\mu(s)}\, ds=+\infty.
\end{equation}
\end{definition}
\begin{remark}
Examples of moduli of continuity satisfying the Osgood condition are $\mu(s)=s$  and $\mu(s)= s\log (e+\frac{1}{s}-1)$.
\end{remark}

We state our main result.

\begin{theorem}\label{main_th}
Let $L$ be the operator
\begin{equation}\label{opL}
L=\partial_t+\sum_{j,k=1}^n \partial_{x_j}(a_{j,k}(t,x)\partial_{x_k}) +\sum_{j=1}^n b_j(t,x)\partial_{x_j}+c(t,x),
\end{equation}
where all the coefficients are supposed to be complex valued, defined in $[0,\,T]\times \R^n$, measurable and bounded.
Let $(a_{j,k}(t,x))_{j,k}$ be  a real symmetric matrix and suppose there exists $\lambda_0\in (0,\,1)$ such that
\begin{equation}\label{elliptic}
\sum_{j,k=1}^n a_{j,k}(t,x)\xi_j\xi_k\geq \lambda_0|\xi|^2,
\end{equation}
for all  $(t,x)\in [0,\,T]\times \R^n$ and for all $\xi\in \R^n$. Under this condition $L$ is a backward parabolic operator.
Let ${\mathcal H}$ be the space of functions such that
\begin{equation}\label{H}
\mathcal H=H^1((0,T), L^2(\R^n_x))\cap L^2 ((0,T), H^2(\R^n_x)).
\end{equation}
Let $\mu$ be a modulus of continuity satisfying the Osgood condition.
Suppose that there exist $\alpha\in (0,\,1)$ and $C>0$ such that, 
\begin{itemize}
\item [i)]
for all $j,k=1,\dots, n$,
\begin{equation}\label{Hol}
a_{j,k}\in C^{0,\alpha}([0,T], L^\infty(\R^n_x))\cap L^\infty ([0,T], Lip(\R^n_x));
\end{equation}
\item [ii)] for all $j,k=1,\dots, n$ and for all $t\in(0,T]$,
\begin{equation}\label{main_cond}
\sup_{s_1,\,s_2\in[t,T],\atop x\in\R^n}\frac{|a_{j,k}(s_1,x)-a_{j,k}(s_2,x)|}{\mu(|s_1-s_2|)}\leq C t^{\alpha-1}.
\end{equation}
\end{itemize}

Then $L$ has the $\mathcal H$-uniqueness property, i.e. if $u\in \mathcal H$, $Lu=0$ in   $[0,T]\times \R^n_x$ and $u(0,x)=0$ in $\R^n_x$, then $u=0$ in $[0,T]\times \R^n_x$.

\end{theorem}


\section{Weight function and Carleman estimate} 
We define
\begin{equation}\label{phi}
\phi(t)= \int_{\frac{1}{t}}^1 \frac{1}{\mu(s)}\, ds.
\end{equation}
The function $\phi:[1,+\infty[\to [0,+\infty[$ is a strictly increasing $C^1$ function and, from Osgood condition, it is bijective. Moreover, for all 
$t\in [1,+\infty[$,
$$
\phi'(t)= \frac{1}{t^2\mu(\frac{1}{t})}.
$$
We remark that $\phi'(1)=1$ and $\phi'$ is decreasing in $[1, +\infty[$, so that $\phi$ is a concave function. We remark also that $\phi^{-1}:[0,+\infty[\to
[1,+\infty[$ and, for all $s\in [0,+\infty[$, 
$$
\phi^{-1}(s)\geq 1+s.
$$
We define 
\begin{equation}\label{psi}
\psi_\gamma(\tau)= \phi^{-1}(\gamma\int_0^{\frac{\tau}{\gamma} }{(T-s)^{\alpha-1}}\, ds),
\end{equation}
where $\tau\in[0,\gamma T]$.
$$
\phi(\psi_\gamma(\tau))=\gamma\int_0^{\frac{\tau}{\gamma} }{(T-s)^{\alpha-1}}\, ds
$$
and 
$$
\phi'(\psi_\gamma(\tau))\psi'_\gamma(\tau)=(T-\frac{\tau}{\gamma})^{\alpha -1}.
$$
Then
$$
\psi'_\gamma(\tau)=(T-\frac{\tau}{\gamma})^{\alpha -1} \cdot (\psi_\gamma(\tau))^2\mu(\frac{1}{\psi_\gamma(\tau)}),
$$
i. e. $\psi_\gamma$ is a solution to the differential equation
$$
u'(\tau)= (T-\frac{\tau}{\gamma})^{\alpha -1} \, u^2(\tau)\, \mu(\frac{1}{u(\tau)}).
$$
Finally we set, for $\tau\in[0, \gamma T]$,
\begin{equation}\label{weight}
\Phi_\gamma(\tau)= \int_0^\tau \psi_\gamma (\sigma)\,d\sigma.
\end{equation}
Remark that, with this definition, $\Phi'(\tau)= \psi_\gamma(\tau)$ and 
\begin{equation}\label{diffeq}
\Phi''_\gamma(\tau)= (T-\frac{\tau}{\gamma})^{\alpha -1}\, (\Phi_\gamma' (\tau))^2 \, \mu(\frac{1}{\Phi_\gamma'(\tau)}).
\end{equation}
In particular, for $t\in (0,\frac{T}{2}]$, 
\begin{equation}\label{weight}
\Phi''_\gamma(\gamma(T-t)) = t^{\alpha-1} \, \Phi_\gamma' (\gamma(T-t))\frac{\mu(\frac{1}{\Phi_\gamma'(\gamma(T-t))})}{\frac{1}{\Phi_\gamma'(\gamma(T-t))}}\geq  t^{\alpha-1}  \geq (\frac{T}{2})^{\alpha-1},
\end{equation}
since $\Phi_\gamma' (\gamma(T-t))= \psi_\gamma(\gamma(T-t))\geq 1$ and $\frac{\mu(s)}{s}\geq 1$ for all $s\in(0,1)$.

\bigskip

We can now state the Carleman estimate.

\begin{theorem}\label{car_est}
In the previous hypotheses there exist $\gamma_0>0$, $C>0$ such that
\begin{equation}\label{carleman}
\begin{array}{ll}
\displaystyle
{\int_0^{\frac{T}{2}} e^{\frac{2}{\gamma} \Phi_\gamma(\gamma(T-t))}\|\partial_tu+\sum_{j,k=1}^n\partial _{x_j}(
a_{j,k}(t,x)\partial_{x_k}u)\|^2_{L^2(\R^n_x)}\, dt}\\[0.4 cm]
\qquad \qquad \displaystyle
{\geq C\gamma^{\frac{1}{2}} \int_0^{\frac{T}{2}} e^{\frac{2}{\gamma} \Phi_\gamma(\gamma(T-t))}
(\|\nabla_x u\|^2_{L^2(\R^n_x)}+\gamma^{\frac{1}{2}}\| u\|^2_{L^2(\R^n_x)})\, dt}
\end{array}
\end{equation}
for all $\gamma>\gamma_0$ and for all $u\in C^\infty_0(\R^{n+1})$ such that $\supp u\subseteq [0,\frac{T}{2}]\times \R^n_x$.

\end{theorem}

The way of obtaining the ${\mathcal H}$-uniqueness from the inequality (\ref{carleman}) is a standard procedure, the details of which can be found in \cite[Par. 3.4]{DSP1}. 





\section{Proof of the Carleman estimate} 

\subsection{Littlewood-Paley decomposition}
We will use the so called Lit\-tle\-wood-Paley theory. We refer to \cite{B}, \cite{Ch}, \cite{M} and \cite{BChD} for the details. 
Let $\psi\in C^{\infty}([0,+\infty[, \R)$ such that $\psi$ is non-increasing and 
$$
\psi(t)=1\quad\text{for}\quad 0\leq t\leq\frac{11}{10}, \qquad \psi(t)=0\quad\text{for}\quad   t\geq\frac{19}{10}.
$$
We set, for $\xi\in\R^n$, 
\begin{equation}\label{defchipsi}
\chi(\xi)=\psi(|\xi|), \qquad \varphi(\xi)= \chi(\xi)-\chi(2\xi).
\end{equation}
Given a tempered distribution $u$, the dyadic blocks are defined by
$$
u_0=\Delta_0 u= \chi(D)u ={\mathcal F}^{-1}(\chi(\xi)\hat u(\xi)), 
$$
$$ u_j= \Delta_j u= \varphi(2^{-j}D)u={\mathcal F}^{-1}(\varphi(2^{-j}\xi)\hat u(\xi))\quad \text{if} \quad j\geq 1,
$$
where we have denoted by ${\mathcal F}^{-1}$ the inverse of the Fourier transform. We introduce also the operator 
$$
S_ku=\sum_{j=0}^{k} \Delta_ju=  {\mathcal F}^{-1}(\chi(2^{-k}\xi)\hat u(\xi)).
$$
We recall some well known facts on Littlewood-Paley demposition.

\begin{proposition}{\rm (\cite[Prop. 3.1]{CM})} \label{char_Hs}
Let $s\in\R$. A temperate distribution $u$ is in $H^s$ if and only if, for all $j\in\N$, $\Delta_ju\in L^2$ and 
$$
\sum_{j=0}^{+\infty} 2^{2js}\|\Delta_ju \|_{L^2}^2<+\infty.
$$

Moreover there exists $C>1$, depending only on $n$ and $s$, such that, for all $u\in H^s$,
\begin{equation}\label{charSob}
\frac{1}{C} \sum_{j=0}^{+\infty} 2^{2js}\|\Delta_ju \|_{L^2}^2\leq \|u\|^2_{H^s}\leq C \sum_{j=0}^{+\infty} 2^{2js}\|\Delta_ju \|_{L^2}^2.
\end{equation}
\end{proposition}

\begin{proposition}{\rm (\cite[Lemma 3.2]{GR}).} \label{char_Hs}
A bounded function $a$ is a Lipschitz-continuous function if and only if
$$
\sup_{k\in \N}\|\nabla (S_ k a)\|_{L^\infty}<+\infty.
$$

Moreover there exists $C>0$, depending only on $n$, such that, for all $a\in Lip$ and for all $k\in \N$,
\begin{equation}\label{estLip}
\|\Delta_k a\|_{L^\infty}\leq C\, 2^{-k}\, \|a\|_{Lip}\qquad\text{and}\qquad \|\nabla(S_k a)\|_{L^\infty}\leq C\,  \|a\|_{Lip},
\end{equation}
where $\|a\|_{Lip}=\|a\|_{L^\infty}+\|\nabla a\|_{L^\infty}$.
\end{proposition}

\subsection{Modified Bony's paraproduct}

\begin{definition}
Let $m\in \N\setminus \{0\}$ and  let $a\in L^\infty$. Let $s\in \R $ and let $u\in H^s$. We define
$$
T^m_a u= S_{m-1}aS_{m+1}u + \sum_{k=m-1}^{+\infty} S_k a \Delta_{k+3} u.
$$
\end{definition}

We recall some known facts on modified Bony's paraproduct.

\begin{proposition} {\rm (\cite[Prop. 5.2.1 and Th. 5.2.8]{M}).} \label{parapr}
Let $m\in \N\setminus \{0\}$ and  let $a\in L^\infty$. Let $s\in \R $. 

Then $T^m_a$ maps $H^s$ into $H^s$ and there exists $C>0$ depending only on $n$, $m$ and $s$, such that, for all $u\in H^s$,
\begin{equation}\label{cont_T}
\| T^m_a u\|_{H^s} \leq C \|a\|_{L^\infty}\, \|u\|_{H^s}.
\end{equation}
Let $m\in \N\setminus \{0\}$ and  let $a\in Lip$. 

Then $a-T^m_a$ maps $L^2$ into $H^1$ and there exists $C'>0$ depending only on $n$, $m$, such that, for all $u\in L^2$,
\begin{equation}\label{cont_a-T}
\| au -T^m_a u\|_{H^{1}} \leq C' \|a\|_{Lip}\, \|u\|_{L^2}.
\end{equation}
\end{proposition}

\begin{proposition}{\rm(\cite[Cor. 3.12]{CM})}\label{pos}
Let $m\in \N\setminus \{0\}$ and  let $a\in Lip$. Suppose that, for all $x\in \R^n$, $a(x)\geq \lambda_0>0$.

Then there exists $m$ depending on $\lambda_0$ and $\|a\|_{Lip}$ such that for all $u\in L^2$,
\begin{equation}\label{pos_a}
{\rm Re}\, \langle T^m_au, u\rangle _{L^2, L^2}\geq \frac{\lambda_0}{2}\|u\|_{L^2}. 
\end{equation}
A similar result remains valid for valued functions when $a$ is replaced by a positive definite matrix $(a_{j,k})_{j,k}$.
\end{proposition}

\begin{proposition}\label{adj} {\rm (\cite[Prop. 3.8 and Prop. 3.11]{CM} and \cite[Prop. 3.8]{DSP2})}
Let $m\in \N\setminus \{0\}$ and  let $a\in Lip$. Let $(T^m_a)^*$ be the adjoint operator of $T^m_a$.

Then there exists $C$ depending only on  $n$ and $m$ such that for all $u\in L^2$,
\begin{equation}\label{adj_a}
\|(T^m_a-(T^m_a)^*)\partial_{x_j}u\|_{L^2}\leq C \|a\|_{Lip} \|u\|_{L^2}. 
\end{equation}
\end{proposition}

We end this subsection with a property which will needed in the proof of the Carleman estimate.

\begin{proposition}\label{comm} {\rm(\cite[Prop. 3.8]{DSP2})}
Let $m\in \N\setminus \{0\}$ and  let $a\in Lip$. Denote by $[\Delta_k, T^m_a]$ the commutator  between $\Delta_k$ and $T^m_a$.

Then there exists $C$ depending only on  $n$ and $m$ such that for all $u\in H^1$,
\begin{equation}\label{comm_T}
(\sum_{h=0}^{+\infty} \|\partial_{x_j}([\Delta_k, T^m_a]\partial_{x_k}u)\|^2_{L^2})^{\frac{1}{2}}\leq C \|a\|_{Lip}\|u\|_{H^1}.
\end{equation}
\end{proposition}


\subsection{Approximated Carleman estimate}

We set
$$
v(t,x)=  e^{\frac{1}{\gamma} \Phi_\gamma(\gamma(T-t))} u(t,x).
$$
The Carleman estimate (\ref{carleman}) becomes: there exist $\gamma_0>0$, $C>0$ such that

\begin{equation}\label{carleman2}
\begin{array}{ll}
\displaystyle
{\int_0^{\frac{T}{2}} \|\partial_tv+\sum_{j,k=1}^n\partial _{x_j}(
a_{j,k}(t,x)\partial_{x_k}v)+ \Phi'_\gamma(\gamma(T-t))v \|^2_{L^2(\R^n_x)}\, dt}\\[0.4 cm]
\qquad \qquad \displaystyle
{\geq C\gamma^{\frac{1}{2}} \int_0^{\frac{T}{2}} 
(\|\nabla_x v\|^2_{L^2(\R^n_x)}+\gamma^{\frac{1}{2}}\| u\|^2_{L^2(\R^n_x)})\, dt,}
\end{array}
\end{equation}
for all $\gamma>\gamma_0$ and for all $v\in C^\infty_0(\R^{n+1})$ such that $\supp u\subseteq [0,\frac{T}{2}]\times \R^n_x$.

First of all, using Proposition \ref{pos}, we fix a value for $m$ in such a way that 
\begin{equation}\label{pos_Ta}
{\rm Re}\, \sum_{j,k}\langle T^m_{a_{j,k}}\partial_{x_k}v, \partial_{x_j}v\rangle_{L^2, L^2}\geq \frac{\lambda_0}{2}\|\nabla v\|_{L^2},
\end{equation}
for all $v\in C^\infty_0(\R^{n+1})$ such that $\supp u\subseteq [0,\frac{T}{2}]\times \R^n_x$. Next we consider 
Proposition \ref{parapr}  and in particular from (\ref{cont_a-T}) we deduce that (\ref{carleman2}) will be a consequence of 
\begin{equation}\label{carleman3}
\begin{array}{ll}
\displaystyle
{\int_0^{\frac{T}{2}} \|\partial_tv+\sum_{j,k=1}^n\partial _{x_j}(
T^m_{a_{j,k}}\partial_{x_k}v)+ \Phi'_\gamma(\gamma(T-t))v \|^2_{L^2(\R^n_x)}\, dt}\\[0.4 cm]
\qquad \qquad \displaystyle
{\geq C\gamma^{\frac{1}{2}} \int_0^{\frac{T}{2}} 
(\|\nabla_x v\|^2_{L^2(\R^n_x)}+\gamma^{\frac{1}{2}}\| u\|^2_{L^2(\R^n_x)})\, dt,}
\end{array}
\end{equation} 
since the difference between (\ref{carleman2}) and (\ref{carleman3}) is absorbed by the right  side part of (\ref{carleman3}) with possibly a different value of $C$ and $\gamma_0$. With a similar argument, using (\ref{charSob}) and (\ref{comm_T}), (\ref{carleman3}) will be deduced from

\begin{equation}\label{carleman4}
\begin{array}{ll}
\displaystyle
{\int_0^{\frac{T}{2}}\sum_{h=0}^{+\infty}\|\partial_tv_h+\sum_{j,k=1}^n\partial _{x_j}(
T^m_{a_{j,k}}\partial_{x_k}v_h)+ \Phi'_\gamma(\gamma(T-t))v_h \|^2_{L^2(\R^n_x)}\, dt}\\[0.4 cm]
\qquad \qquad \displaystyle
{\geq C\gamma^{\frac{1}{2}} \int_0^{\frac{T}{2}} 
\sum_{h=0}^{+\infty}(\|\nabla_x v_h\|^2_{L^2(\R^n_x)}+\gamma^{\frac{1}{2}}\| v_h\|^2_{L^2(\R^n_x)})\, dt,}
\end{array}
\end{equation} 
where we have denoted by $v_h$ the dyadic block $\Delta_h v$.

We fix our attention on each of the terms
$$
\int_0^{\frac{T}{2}} \|\partial_tv_h+\sum_{j,k=1}^n\partial _{x_j}(
T^m_{a_{j,k}}\partial_{x_k}v_h)+ \Phi'_\gamma(\gamma(T-t))v_h \|^2_{L^2(\R^n_x)}\, dt.
$$
We have
\begin{equation}\label{carl_v_h}
\begin{array}{ll}
\displaystyle{\int_0^{\frac{T}{2}} \|\partial_tv_h+\sum_{j,k=1}^n\partial _{x_j}(
T^m_{a_{j,k}}\partial_{x_k}v_h)+ \Phi'_\gamma(\gamma(T-t))v_h \|^2_{L^2(\R^n_x)}\, dt}\\[0.4 cm]
\displaystyle{= \int_0^{\frac{T}{2}} \Big( \|\partial_tv_h\|^2_{L^2}+\|\sum_{j,k=1}^n\partial _{x_j}(
T^m_{a_{j,k}}\partial_{x_k}v_h) + \Phi'_\gamma(\gamma(T-t))v_h \|^2_{L^2}}\\[0.4 cm]
\quad\displaystyle{+ \gamma\Phi''_\gamma(\gamma(T-t))\|v_h \|^2_{L^2}
+2\,{\rm Re}\, \langle \partial_tv_h, \sum_{j,k=1}^n\partial _{x_j}(
T^m_{a_{j,k}}\partial_{x_k}v_h)\rangle _{L^2, L^2} \Big)\, dt}
\end{array}
\end{equation}
Let consider the last term in (\ref{carl_v_h}). We  define, for $\varepsilon\in (0, \frac{T}{2}]$, 
$$
{\tilde a}_{j,k,  \varepsilon}(t,x)=\left\{
\begin{array}{ll}
a_{j,k} (t,x),\qquad& \text{if}\ \ t\geq \varepsilon,\\[0.2cm]
a_{j,k} (\varepsilon,x),\qquad& \text{if}\ \ t<\varepsilon,
\end{array}
\right.
$$
and
$$
a_{j,k,\varepsilon}(t,x)= \int_{-\varepsilon}^{\varepsilon} \rho_\varepsilon(s) {\tilde a}_{j,k,\varepsilon}(t-s, x)\,ds,
$$
where $\rho\in C^\infty_0(\R)$ with $\supp \rho\subseteq [-1,\,1]$, $\int_\R \rho(s)ds=1$, $\rho (s)\geq 0$ and $\rho_\varepsilon(s)= \frac{1}{\varepsilon}
\rho (\frac{s}{\varepsilon})$. With a  strightforward computation, form (\ref{Hol}) and (\ref{main_cond}), we obtain
\begin{equation}\label{a-a_e}
|a_{j,k}(t,x)-a_{j,k, \varepsilon}(t,x)|\leq C\, \min\{ \varepsilon^\alpha,\, t^{\alpha-1}\, \mu(\varepsilon)\},
\end{equation}
and
\begin{equation}\label{a'_e}
|\partial_t a_{j,k, \varepsilon}(t,x)|\leq C\min\{ \varepsilon^{\alpha-1},\, t^{\alpha-1}\,\frac{ \mu(\varepsilon)}{\varepsilon }\},
\end{equation}
for all $j,k= \,\dots, n$ and for all $(t,x)\in [0, \frac{T}{2}]\times \R^n_x$. We deduce
$$
\begin{array}{ll}
\displaystyle{\int_0^{\frac{T}{2}}2\,{\rm Re}\, \langle \partial_tv_h, \sum_{j,k=1}^n\partial _{x_j}(
T^m_{a_{j,k}}\partial_{x_k}v_h)\rangle _{L^2, L^2} \, dt}\\[0.4cm]
\displaystyle{=-2\,{\rm Re}\, \int_0^{\frac{T}{2}} \sum_{j,k=1}^n \langle \partial_{x_j} \partial_tv_h,\;
T^m_{a_{j,k}} \partial _{x_k} v_h\rangle _{L^2, L^2} \, dt}\\[0.4cm]
\displaystyle{=-2\,{\rm Re}\, \int_0^{\frac{T}{2}} \sum_{j,k=1}^n \langle \partial_{x_j} \partial_tv_h,\;
(T^m_{a_{j,k}}-T^m_{a_{j,k,\varepsilon }}) \partial _{x_k} v_h\rangle _{L^2, L^2} \, dt}\\[0.4cm]
\qquad\qquad\qquad \displaystyle{-2\,{\rm Re}\, \int_0^{\frac{T}{2}} \sum_{j,k=1}^n \langle \partial_{x_j} \partial_tv_h,\;
T^m_{a_{j,k,\varepsilon }} \partial _{x_k} v_h\rangle _{L^2, L^2} \, dt}.
\end{array}
$$

Now, $T^m_{a_{j,k}}-T^m_{a_{j,k,\varepsilon }}= T^m_{a_{j,k}-a_{j,k,\varepsilon }}$ and, from (\ref{cont_T}) and (\ref{a-a_e}), 
$$
\begin{array}{ll}
\displaystyle{
\|(T^m_{a_{j,k}}-T^m_{a_{j,k,\varepsilon }}) \partial _{x_k} v_h\|_{L^2}}&\displaystyle{=\|T^m_{a_{j,k}-a_{j,k,\varepsilon }}\partial _{x_k} v_h\|_{L^2}}\\[0.4cm]
&\displaystyle{\leq C\, \min\{\varepsilon^\alpha, t^{\alpha-1}\mu(\varepsilon)\} \|\partial _{x_k} v_h\|_{L^2}.}
\end{array}
$$
Moreover $\|\partial_{x_j}v_h\|_{L^2}\leq 2^{h+1} \|v_h\|_{L^2}$ and $\|\partial_{x_j}\partial_t v_h\|_{L^2}\leq 2^{h+1} \| \partial_t  v_h\|_{L^2}$, so that
$$
\begin{array}{ll}
\displaystyle{| 2\,{\rm Re}\, \int_0^{\frac{T}{2}} \sum_{j,k=1}^n \langle \partial_{x_j} \partial_tv_h,\;
(T^m_{a_{j,k}}-T^m_{a_{j,k,\varepsilon }}) \partial _{x_k} v_h\rangle _{L^2, L^2} \, dt}\\[0.4cm]
\quad\displaystyle{\leq 2C \, \int_0^{\frac{T}{2}} \min\{\varepsilon^\alpha, t^{\alpha-1}\mu(\varepsilon)\} \sum_{j,k=1}^n
\|\partial_{x_j}\partial_t v_h\|_{L^2}\|\partial_{x_k}v_h\|_{L^2}\,dt}\\[0.4cm]
\quad\displaystyle{\leq \frac{C}{N}  \, \int_0^{\frac{T}{2}} \| \partial_t  v_h\|^2_{L^2}\,dt
+CN\, 2^{4(h+1)} \, \int_0^{\frac{T}{2}}  \min\{\varepsilon^\alpha, t^{\alpha-1}\mu(\varepsilon)\} \|v_h\|^2_{L^2}\, dt},
\end{array}
$$
where $C$ depends only on $n$, $m$ and $\|a_{j,k}\|_{L^\infty}$ and $N>0$ can be chosen arbitrarily. 

Similarly
$$
\begin{array}{ll}
\displaystyle{-2\,{\rm Re}\, \int_0^{\frac{T}{2}} \sum_{j,k=1}^n \langle \partial_{x_j} \partial_tv_h,\;
T^m_{a_{j,k,\varepsilon }} \partial _{x_k} v_h\rangle _{L^2, L^2} \, dt}\\[0.4 cm]
\qquad\quad\displaystyle{=\int_0^{\frac{T}{2}} \sum_{j,k=1}^n \langle \partial_{x_j} v_h,\; T^m_{ \partial_t a_{j,k,\varepsilon }} \partial _{x_k} v_h\rangle _{L^2, L^2} \, dt}\\[0.4 cm]
\qquad\qquad\qquad\displaystyle{+ \int_0^{\frac{T}{2}} \sum_{j,k=1}^n \langle \partial_{x_j} v_h,\; (T^m_{a_{j,k,\varepsilon }}-(T^m_{a_{j,k,\varepsilon }})^*) \partial _{x_k} \partial_t v_h\rangle _{L^2, L^2} \, dt}.
\end{array}
$$
From (\ref{cont_T}) and (\ref{a'_e}) we have
$$
\begin{array}{ll}
\displaystyle{|\int_0^{\frac{T}{2}} \sum_{j,k=1}^n \langle \partial_{x_j} v_h,\; T^m_{ \partial_t a_{j,k,\varepsilon }} \partial _{x_k} v_h\rangle _{L^2, L^2} \, dt|}\\[0.4cm]
\qquad\qquad
\displaystyle{ \leq C\, 2^{2(h+1)} \int_0^{\frac{T}{2}} \min\{\varepsilon^{\alpha-1}, t^{\alpha-1}\frac{\mu(\varepsilon)}{\varepsilon}\} \|v_h\|^2_{L^2}\, dt},
\end{array}
$$
and, from (\ref{adj_a}),
$$
\begin{array}{ll}
\displaystyle{| \int_0^{\frac{T}{2}} \sum_{j,k=1}^n \langle \partial_{x_j} v_h,\; (T^m_{a_{j,k,\varepsilon }}-(T^m_{a_{j,k,\varepsilon }})^*) \partial _{x_k} \partial_t v_h\rangle _{L^2, L^2} \, dt|}\\[0.4 cm]
\qquad\displaystyle{\leq C  \int_0^{\frac{T}{2}}  \|\nabla v_h\|_{L^2} \| \partial_t v_h\|_{L^2}\, dt}\\[0.4 cm]
\qquad\qquad \displaystyle{\leq  \frac{C}{N}  \int_0^{\frac{T}{2}}  \| \partial_t v_h\|_{L^2}^2\,dt + CN  \, 2^{2(h+1)} \int_0^{\frac{T}{2}}  \| v_h\|_{L^2}^2\,dt},
\end{array}
$$
where $C$ depends only on $n$, $m$ and $\|a_{j,k}\|_{Lip}$ and $N>0$ can be chosen arbitrarily. 

As a conclusion, from (\ref{carl_v_h}), we finally obtain
\begin{equation}\label{carl_v_h_2}
\begin{array}{ll}
\displaystyle{\int_0^{\frac{T}{2}} \|\partial_tv_h+\sum_{j,k=1}^n\partial _{x_j}(
T^m_{a_{j,k}}\partial_{x_k}v_h)+ \Phi'_\gamma(\gamma(T-t))v_h \|^2_{L^2(\R^n_x)}\, dt}\\[0.4 cm]
\displaystyle{\geq \int_0^{\frac{T}{2}} \Big( \| \sum_{j,k=1}^n\partial _{x_j}(
T^m_{a_{j,k}}\partial_{x_k}v_h) + \Phi'_\gamma(\gamma(T-t))v_h \|^2_{L^2} +\gamma \Phi''_\gamma(\gamma(T-t))\|v_h \|^2_{L^2} }\\[0.4 cm]
-C ( 2^{4(h+1)} \min\{\varepsilon^{\alpha}, t^{\alpha-1}\mu(\varepsilon)\} +
2^{2(h+1)}(\min\{\varepsilon^{\alpha-1}, t^{\alpha-1}\frac{\mu(\varepsilon)}{\varepsilon}\} + 1)) |v_h \|^2_{L^2}\Big)\, dt.\\[0.4 cm]
\end{array}
\end{equation}
\subsection{End of the proof}

We start considering (\ref{carl_v_h_2}) for $h=0$. We fix $\varepsilon=\frac{1}{2}$. Recalling (\ref{weight}) we have
\begin{equation*}
\begin{array}{ll}
\displaystyle{\int_0^{\frac{T}{2}} \|\partial_tv_h+\sum_{j,k=1}^n\partial _{x_j}(
T^m_{a_{j,k}}\partial_{x_k}v_0)+ \Phi'_\gamma(\gamma(T-t))v_0 \|^2_{L^2(\R^n_x)}\, dt}\\[0.4 cm]
\qquad\qquad \displaystyle{\geq \int_0^{\frac{T}{2}} (\gamma \Phi''_\gamma(\gamma(T-t))- C')\|v_0 \|^2_{L^2} }\\[0.4 cm]
\qquad\qquad\qquad\qquad \displaystyle{\geq \int_0^{\frac{T}{2}} (\gamma (\frac{T}{2})^{\alpha-1}- C')\|v_0 \|^2_{L^2} \, dt}.
\end{array}
\end{equation*}
Choosing a suitable $\gamma_0$, we have that, for all $\gamma>\gamma_0$,
\begin{equation}\label{carl_v_0}
\int_0^{\frac{T}{2}} \|\partial_tv_h+\sum_{j,k=1}^n\partial _{x_j}(
T^m_{a_{j,k}}\partial_{x_k}v_0)+ \Phi'_\gamma(\gamma(T-t))v_0 \|^2_{L^2}\, dt \geq \frac{\gamma}{2} \int_0^{\frac{T}{2}} 
\|v_0 \|^2_{L^2}\, dt.
\end{equation}

We consider (\ref{carl_v_h_2}) for $h\geq 1$. We fix $\varepsilon=2^{-2h}$. We have
$$
\begin{array}{ll}
\displaystyle{\int_0^{\frac{T}{2}} \|\partial_tv_h+\sum_{j,k=1}^n\partial _{x_j}(
T^m_{a_{j,k}}\partial_{x_k}v_h)+ \Phi'_\gamma(\gamma(T-t))v_h \|^2_{L^2(\R^n_x)}\, dt}\\[0.4 cm]
\displaystyle{\geq \int_0^{\frac{T}{2}} \Big( \| \sum_{j,k=1}^n\partial _{x_j}(
T^m_{a_{j,k}}\partial_{x_k}v_h) + \Phi'_\gamma(\gamma(T-t))v_h \|^2_{L^2}  }\\[0.4 cm]
\quad \displaystyle{+\big(\gamma \Phi''_\gamma(\gamma(T-t))-C ( 2^{4h} \min\{2^{-2h\alpha}, t^{\alpha-1}\mu(2^{-2h})\} +2^{2h})\big)
\|v_h \|^2_{L^2}\Big)\, dt}\\[0.4 cm]
\displaystyle{\geq \int_0^{\frac{T}{2}} \Big( (\| \sum_{j,k=1}^n\partial _{x_j}(
T^m_{a_{j,k}}\partial_{x_k}v_h)\|_{L^2} - \Phi'_\gamma(\gamma(T-t))\|v_h \|^2_{L^2} \big)^2 }\\[0.4 cm]
\quad \displaystyle{+\big(\gamma \Phi''_\gamma(\gamma(T-t))-C ( 2^{4h} \min\{2^{-2h\alpha}, t^{\alpha-1}\mu(2^{-2h})\} +2^{2h})\big)
\|v_h \|^2_{L^2}\Big)\, dt.}\\[0.4 cm]
\end{array}
$$
From (\ref{pos_Ta}) it is possible to deduce that 
\begin{equation}\label{Ta_v_h}
\| \sum_{j,k=1}^n\partial _{x_j}(
T^m_{a_{j,k}}\partial_{x_k}v_h)\|_{L^2}\geq \frac{\lambda_0}{8}\, 2^{2h} \|v_h\|^2_{L^2}.
\end{equation}

Suppose first that 
$$
\Phi'_\gamma (\gamma(T-t))\leq \frac{\lambda_0}{16} 2^{2h}.
$$
From (\ref{pos_Ta}) we have
$$
\| \sum_{j,k=1}^n\partial _{x_j}(
T^m_{a_{j,k}}\partial_{x_k}v_h)\|_{L^2} - \Phi'_\gamma(\gamma(T-t))\|v_h \|^2_{L^2}
\geq \frac{\lambda_0}{16}\, 2^{2h} \|v_h \|^2_{L^2}
$$
and then, using also (\ref{weight}), we obtain
$$
\begin{array}{ll}
\displaystyle{\int_0^{\frac{T}{2}} \|\partial_tv_h+\sum_{j,k=1}^n\partial _{x_j}(
T^m_{a_{j,k}}\partial_{x_k}v_h)+ \Phi'_\gamma(\gamma(T-t))v_h \|^2_{L^2(\R^n_x)}\, dt}\\[0.4 cm]
\displaystyle{\geq \int_0^{\frac{T}{2}} \Big( (\| \sum_{j,k=1}^n\partial _{x_j}(
T^m_{a_{j,k}}\partial_{x_k}v_h)\|_{L^2} - \Phi'_\gamma(\gamma(T-t))\|v_h \|^2_{L^2} \big)^2 }\\[0.4 cm]
\quad \displaystyle{+\big(\gamma \Phi''_\gamma(\gamma(T-t))-C ( 2^{4h} \min\{2^{-2h\alpha}, t^{\alpha-1}\mu(2^{-2h})\} +2^{2h})\big)
\|v_h \|^2_{L^2}\Big)\, dt}\\[0.4 cm]
\displaystyle{\geq \int_0^{\frac{T}{2}} \Big( (\frac{\lambda_0}{16}\, 2^{2h} )^2 +\gamma (\frac{T}{2})^{\alpha-1}-C ( 2^{(4-2\alpha) h} \big)
\|v_h \|^2_{L^2}\Big)\, dt.}
\end{array}
$$
Then there exist  $\gamma_0>0$ and  $C>0$ such that, for all $\gamma\geq \gamma_0$ and for all $h\geq 1$,
\begin{equation}\label{final1}
\begin{array}{ll}
\displaystyle{\int_0^{\frac{T}{2}} \|\partial_tv_h+\sum_{j,k=1}^n\partial _{x_j}(
T^m_{a_{j,k}}\partial_{x_k}v_h)+ \Phi'_\gamma(\gamma(T-t))v_h \|^2_{L^2(\R^n_x)}\, dt}\\[0.4 cm]
\qquad\qquad\qquad\qquad \displaystyle{\geq C \int_0^{\frac{T}{2}}
(\gamma + \gamma^{\frac{1}{2}}2^{2h}) \|v_h \|^2_{L^2}\Big)\, dt}
\end{array}
\end{equation}

Suppose finally that 
$$
\Phi'_\gamma (\gamma(T-t))\geq \frac{\lambda_0}{16} 2^{2h}.
$$
From (\ref{diffeq}), the fact that $\lambda_0 \leq 1$ and the properties of the modulus of continuity $\mu$
$$
\begin{array}{ll}
\displaystyle{\Phi''(\gamma(T-t))}&\displaystyle{=t^{\alpha-1} (\Phi_\gamma' (\gamma(T-t)))^2 \, \mu(\frac{1}{\Phi_\gamma'(\gamma(T-t)))})}\\[0.4 cm]
&\displaystyle{\geq t^{\alpha-1} (\frac{\lambda_0}{16})^2  2^{4h}\mu(\frac {16}{\lambda_0} 2^{-2h}) 
\geq t^{\alpha-1} (\frac{\lambda_0}{16})^2  2^{4h}\mu(2^{-2h}).}
\end{array}
$$
and 
$$
\begin{array}{ll}
\displaystyle{\Phi''(\gamma(T-t))}&\displaystyle{=t^{\alpha-1} (\Phi_\gamma' (\gamma(T-t)))^2 \, \mu(\frac{1}{\Phi_\gamma'(\gamma(T-t)))})}\\[0.4 cm]
&\displaystyle{= t^{\alpha-1} \, \Phi_\gamma' (\gamma(T-t))\frac{\mu(\frac{1}{\Phi_\gamma'(\gamma(T-t))})}{\frac{1}{\Phi_\gamma'(\gamma(T-t))}}\geq (\frac{T}{2})^{\alpha-1}.}
\end{array}
$$
Consequently
$$
\begin{array}{ll}
\displaystyle{\int_0^{\frac{T}{2}} \|\partial_tv_h+\sum_{j,k=1}^n\partial _{x_j}(
T^m_{a_{j,k}}\partial_{x_k}v_h)+ \Phi'_\gamma(\gamma(T-t))v_h \|^2_{L^2(\R^n_x)}\, dt}\\[0.4 cm]
 \displaystyle{\geq\int_0^{\frac{T}{2}} \big(\gamma \Phi''_\gamma(\gamma(T-t))-C ( 2^{4h} \min\{2^{-2h\alpha}, t^{\alpha-1}\mu(2^{-2h})\} +2^{2h})\big)
\|v_h \|^2_{L^2}\, dt}\\[0.4 cm]
\displaystyle{\geq \int_0^{\frac{T}{2}} \frac{\gamma}{2} (t^{\alpha-1} \big(\frac{\lambda_0}{16})^2  2^{4h}\mu(2^{-2h})+
(\frac{T}{2})^{\alpha-1}\big)
-C \big( t^{\alpha-1}2^{4 h}\mu(2^{-2h}) +2^{2h} \big)
\|v_h \|^2_{L^2}\Big)\, dt.}
\end{array}
$$
Then there exist  $\gamma_0>0$ and  $C>0$ such that, for all $\gamma\geq \gamma_0$ and for all $h\geq 1$,
\begin{equation}\label{final2}
\begin{array}{ll}
\displaystyle{\int_0^{\frac{T}{2}} \|\partial_tv_h+\sum_{j,k=1}^n\partial _{x_j}(
T^m_{a_{j,k}}\partial_{x_k}v_h)+ \Phi'_\gamma(\gamma(T-t))v_h \|^2_{L^2(\R^n_x)}\, dt}\\[0.4 cm]
\qquad\qquad\qquad\qquad \displaystyle{\geq C \gamma \int_0^{\frac{T}{2}}
(1 + 2^{2h}) \|v_h \|^2_{L^2}\Big)\, dt.}
\end{array}
\end{equation}
As a conclusion, form (\ref{carl_v_0}), (\ref{final1}) and (\ref{final2}), there exist  $\gamma_0>0$ and  $C>0$ such that, for all $\gamma\geq \gamma_0$ and for all $h\in \N$,
\begin{equation}\label{final3}
\begin{array}{ll}
\displaystyle{\int_0^{\frac{T}{2}} \|\partial_tv_h+\sum_{j,k=1}^n\partial _{x_j}(
T^m_{a_{j,k}}\partial_{x_k}v_h)+ \Phi'_\gamma(\gamma(T-t))v_h \|^2_{L^2(\R^n_x)}\, dt}\\[0.4 cm]
\qquad\qquad\qquad\qquad \displaystyle{\geq C \int_0^{\frac{T}{2}}
(\gamma + \gamma^{\frac{1}{2}}2^{2h}) \|v_h \|^2_{L^2}\Big)\, dt}
\end{array}
\end{equation}
and (\ref{carleman4}) follows. The proof is complete.





\section{A counterexample}
\begin{theorem}
There exists 
$$
l\in \big(\bigcap_{\alpha\in [0,1[}  C^{0,\alpha}(\R)\big)\cap C^\infty(\R\setminus\{0\})
$$
with
\begin{equation}\label{ellip}
\frac{1}{2}\leq l(t)\leq \frac{3}{2},\qquad \text{for all}\ \  t\in \R,
\end{equation}

\begin{equation}\label{osc}
|l'(t)|\leq C_\varepsilon |t|^{-(1+\varepsilon)},\qquad \text{for all}\ \  \varepsilon >0 \ \ \text{and}\ \ t\in \setminus\{0\},
\end{equation}

\smallskip
\noindent
and there exist $u, \; b_1, \; b_2, \; c\in C^\infty_b(\R_t\times \R^2_x)$, with 
$$\supp u=\{(t,x)\in \R_t\times \R^2_x\,\big| \, t\geq 0\},$$ such that
\begin{equation*}
\partial_t u +\partial^2_{x_1} u+ l\partial^2_{x_2}u+b_1\partial_{x_1}u+b_2\partial_{x_2}u+cu=0 \qquad \text{in}\ \ \R_t\times \R^2_x .
\end{equation*}
\end{theorem}

\begin{remark}
Actually the function $l$ will satisfy
\begin{equation}\label{osc1}
\sup_{t\not =0} (\frac{|t|}{1+|\log|t||})|l'(t)|<+\infty.
\end{equation}
From (\ref{osc1}) it is easy to obtain (\ref{osc}).
\end{remark}
\begin{proof}
We will follow the proof of Theorem 1 in \cite{P} (see also Theorem 3 in \cite {DSP1}). Let $A, \; B, \; C,\; J$ be four $C^\infty$ functions, defined in $\R$, with 
$$
0\le A(s),\, B(s),\; C(s)\le 1 \quad \text{and}\quad -2\le J(s)\le 2,\quad \text {for all} \ \  s\in\R,
$$
and
$$
\begin{array}{ll}
A(s)=1,\quad \text{for}\ \ s\le\frac{1}{5},\qquad\qquad& A(s)=0,\quad \text{for}\ \ s\ge\frac{1}{4},\\[0.2 cm]
B(s)=0,\quad \text{for}\ \ s\le0 \ \ \text{or}\ \ s\ge 1,\qquad\qquad& B(s)=1,\quad \text{for}\ \ \frac{1}{6}\leq s\le\frac{1}{2},\\[0.2 cm]
C(s)=0,\quad \text{for}\ \ s\le\frac{1}{4},\qquad\qquad& C(s)=1,\quad \text{for}\ \ s\ge\frac{1}{3},\\[0.2 cm]
J(s)=-2,\quad \text{for}\ \ s\le\frac{1}{6} \ \ \text{or}\ \ s\ge \frac{1}{2},\qquad\qquad& J(s)=2,\quad \text{for}\ \ \frac{1}{5}\leq s\le\frac{1}{3}.
\end{array}
$$
Let $(a_n)_n,\; (z_n)_n$ be two real sequences such that
\begin{equation}\label{seq a}
-1<a_n<a_{n+1}, \quad\text{for all} \ \ n\geq 1, \quad\text{and}\quad \lim_n a_n=0,
\end{equation}
\begin{equation}\label{seq z}
1<z_n<z_{n+1}, \quad\text{for all} \ \ n\geq 1, \quad\text{and}\quad \lim_n z_n=+\infty.
\end{equation}
We define
$$
\begin{array}{c}
r_n= a_{n+1}-a_n,\\[0.4 cm]
q_1=0 \quad\text{and}\quad \displaystyle{q_n=\sum_{k=2}^n z_kr_{k-1},} \quad\text{for}\ \ n\geq 2,\\[0.4 cm]
p_n=(z_{n+1}-z_n)r_n.
\end{array}
$$
We require
\begin{equation}\label{seq p}
p_n>1,\quad\text{for all}\ \ n\geq 1.
\end{equation}
We set
$$
A_n(t)= A(\frac{t-a_n}{r_n}),\quad B _n(t)= B(\frac{t-a_n}{r_n}),
$$
$$
C_n(t)= C(\frac{t-a_n}{r_n}),\quad J_n(t)= J(\frac{t-a_n}{r_n}).
$$

We define
$$
\begin{array}{ll}
v_n(t,x_1)= \exp(-q_n-z_n(t-a_n))\cos \sqrt{z_n} \,x_1,\\[0.2 cm]
w_n(t,x_2)= \exp(-q_n-z_n(t-a_n)+J_n(t)p_n)\cos \sqrt{z_n} \,x_2,
\end{array}
$$

$$
\begin{array}{ll}
u(t, x_1,x_2)\\[0.4cm]
=\left\{
\begin{array}{ll}
v_1(t, x_1),&\text{for} \ \ t\leq a_1,\\[0.4 cm]
A_n(t)v_n(t,x_1)+B_n(t)w_n(t,x_2)+C_n(t)v_{n+1}(t,x_1),  & \text{for} \ \ a_n\leq t\leq a_{n+1},\\[0.4 cm]
0,& \text{for} \ \ t\geq 0.\end{array}
\right.
\end{array}
$$
The condition 
\begin{equation}\label{cond1}
\lim_n \, \exp(-q_n+2p_n) z_{n+1}^\alpha p_n^\beta r_n^{-\gamma}=0, \qquad\text{for all}\ \ 	\alpha,\,\beta,\,\gamma>0,
\end{equation}
implies that $u\in C^\infty_b(\R_t\times \R^2_x)$.

We define 
$$
l(t)=\left\{
\begin{array}{ll}
1,&\text{for} \ \ t\leq a_1\ \ \text{or}\ \ t\geq 0,\\[0.4 cm]
1+J'_n(t)p_nz_n^{-1},\qquad& \text{for}  \ \ a_n\leq t\leq a_{n+1}.
\end{array}
\right.
$$
$l$ is a $C^\infty(\R\setminus\{0\})$ function. The condition 
\begin{equation}\label{cond2}
\sup_n \{p_nr_n^{-1} z_n^{-1}\}\leq \frac{1}{2\|J'\|_{L^\infty}}
\end{equation}
implies (\ref{ellip}), i. e. the operator 
$$
L=\partial_t-\partial^2_{x_1}-l(t)\partial^2_{x_2}
$$
is a parabolic operator. Moreover  $l$ is in  $\bigcap_{\alpha\in [0,1[} C^{0,\alpha}(\R)$ if
\begin{equation}\label{cond3}
\sup_n\{p_nr_n^{-1-\alpha}z_n^{-1} \}<+\infty, \qquad \text{for all}\ \ \alpha\in [0,1[.
\end{equation}

Finally, we define
$$
\begin{array}{l}
\displaystyle{b_1=-\frac{Lu}{u^2+(\partial_{x_1}u)^2+(\partial_{x_2}u)^2}\partial_{x_1}u,}\\[0.4 cm]
\displaystyle{b_2=-\frac{Lu}{u^2+(\partial_{x_1}u)^2+(\partial_{x_2}u)^2}\partial_{x_2}u,}\\[0.4 cm]
\displaystyle{c=-\frac{Lu}{u^2+(\partial_{x_1}u)^2+(\partial_{x_2}u)^2}u.}
\end{array}
$$
As in \cite{P} and \cite{DSP1}, the functions $b_1,\, b_2, \, c$ are in $C^\infty_b(\R_t\times \R^2_x)$ if
\begin{equation}\label{cond4}
\lim_n \, \exp(-p_n) z_{n+1}^\alpha p_n^\beta r_n^{-\gamma}=0, \qquad\text{for all}\ \ 	\alpha,\,\beta,\,\gamma>0.
\end{equation}

We choose, for  $j_0\geq 2$,  
$$
a_n=-{e^{-\sqrt{\log (n+j_0)}}},\qquad\qquad z_n=(n+j_0)^3.
$$
With this choice (\ref{seq a}) and (\ref{seq z}) are satisfied and  we have
\begin{equation*}
r_n\sim e^{-\sqrt{\log (n+j_0)}}\, \frac{1}{(n+j_0)\sqrt{\log (n+j_0)}},
\end{equation*}
where, for sequences $(f_n)_n,\, (g_n)_n$,  $f_n\sim g_n$ means $\lim_n \frac{f_n}{g_n} = \lambda$, for some $\lambda>0$.
Similarly
\begin{equation*}
p_n\sim e^{-\sqrt{\log (n+j_0)}}\, \frac{n+j_0}{\sqrt{\log (n+j_0)}}
\end{equation*}
and condition (\ref{seq p})  is verified, for a suitable fixed $j_0$. Remarking that 
we have, for $j_0$ suitably large, 
$$
q_n = \sum_{k=2}^n z_kr_{k-1}\geq z_nr_{n-1}\geq \lambda (n+j_0)^{\frac{7}{4}}
$$
and 
$$
p_n\leq \lambda (n+j_0)^{\frac{5}{4}}
$$
for some $\lambda>0$. Finally
$$
p_n r_n^{-1} z_n^{-1}\sim  \frac{1}{n+j_0}.
$$
As a consequence (\ref{cond1}), (\ref{cond2}), (\ref{cond3}) and (\ref{cond4}) are satisfied for a suitable fixed $j_0$. It remains to check (\ref{osc1}).
We have
$$
|l'(t)|\leq \|J''\|_{L^\infty} p_n\, r_n^{-2}\,z_n^{-1},\qquad \text{for}  \ \ a_n\leq t\leq a_{n+1}
$$
and consequently
$$
\begin{array}{ll}
\displaystyle{\sup_{t\not =0} (\frac{|t|}{1+|\log|t||})|l'(t)| }&=\  \displaystyle{ \sup_n \sup_{t\in [a_n,a_{n+1}]} (\frac{|t|}{1+|\log|t||})|l'(t)|}\\[0.8 cm]
&\leq\  \displaystyle{\sup_n \, (\frac{a_n}{1-\log a_n}) \|J''\|_{L^\infty} p_n\, r_n^{-2}\,z_n^{-1}}\\[0.8 cm]
&\leq C.
\end{array}
$$
The conclusion of the theorem is reached simply exchanging $t$ with $-t$.
\end{proof}



\begin{thebibliography}{99}

\bibitem{BChD} Bahouri, Hajer; Chemin, Jean-Yves; Danchin, Rapha\'el \lq\lq Fourier analysis and nonlinear partial differential equations\rq\rq. Grundlehren der Mathe\-ma\-ti\-schen Wissenschaften, 343. Springer, Heidelberg, 2011. 
  
\bibitem{B} Bony, Jean-Michel  {\it  Calcul symbolique et propagation des singularit\'es pour les \'equations aux d\'eriv\'ees partielles non lin\'eaires}.  Ann. Sci. \'Ecole Norm. Sup. (4) 14 (1981), no. 2, 209--246.    

\bibitem{Ch} Chemin, Jean-Yves \lq\lq Fluides parfaits incompressibles\rq\rq.  Ast\'erisque, 230. Soci\'et\'e Math\'ematique de France, Paris, 1995.





\bibitem{CM}  Colombini, Ferruccio; M\'etivier, Guy {\it  The Cauchy problem for wave equations with non Lipschitz coefficients; application to continuation of solutions of some nonlinear wave equations}. Ann. Sci. \'Ec. Norm. Sup\'er. (4) 41 (2008), no. 2, 177--220.

\bibitem {DSP1} Del Santo, Daniele; Prizzi, Martino {\it Backward uniqueness for parabolic operators whose coefficients are non-Lipschitz continuous in time}. J. Math. Pures Appl. (9) 84 (2005), no. 4, 471--491.
 
\bibitem {DSP2} Del Santo, Daniele; Prizzi, Martino {\it A new result on backward uniqueness for parabolic operators}. Ann. Mat. Pura Appl. (4) 194 (2015), no. 2, 387--403.

\bibitem{DJ}  Del Santo, Daniele; J\"ah,  Christian {\it Non-uniqueness  and  uniqueness  in  the  Cauchy  problem  of  elliptic and backward-parabolic equations}, in "Progress in Partial Differential Equations -  Asymptotic  Profiles,  Regularity  and  Well-Posedness",  M.  Ruzhansky,  M.  Reissig eds.,  Springer  Proceedings  in  Mathematics and  Statistics 44,  Springer  International Publishing, Basel 2013, pp.  27--52.

\bibitem{GR}  G\'erard, Patrick; Rauch, Jeffrey {\it Propagation de la r\'egularit\'e locale de solutions d'\'equations hyperboliques non lin\'eaires}. (French) [Propagation of the local regularity of solutions of nonlinear hyperbolic equations] Ann. Inst. Fourier (Grenoble) 37 (1987), no. 3, 65--84.


\bibitem{M} M\'etivier, Guy \lq\lq Para-differential calculus and applications to the Cauchy problem for nonlinear systems\rq\rq. Centro di Ricerca Matematica Ennio De Giorgi (CRM) Series, 5. Edizioni della Normale, Pisa, 2008.


\bibitem {P}  Pli\'s, Andrzej {\it On non-uniqueness in Cauchy problem for an elliptic second order differential equation}. Bull. Acad. Polon. Sci. S\`er. Sci. Math. Astronom. Phys. 11 (1963), 95--100.


\end{thebibliography}
\end{document}